\documentclass[11pt]{article}

\usepackage[english]{babel}
\usepackage[utf8x]{inputenc}
\usepackage[T1]{fontenc}

\usepackage[a4paper,top=2.54cm,bottom=2.54cm,left=2.54cm,right=2.54cm,marginparwidth=1.75cm]{geometry} 

\usepackage{amsmath}
\usepackage{amsfonts}
\usepackage{amssymb}
\usepackage{graphicx}
\usepackage[colorinlistoftodos]{todonotes}
\usepackage[colorlinks=true, allcolors=blue]{hyperref}
\usepackage{graphicx} 
\usepackage{booktabs} 

\usepackage{caption}
\usepackage{subcaption}
\usepackage{float}
\usepackage{setspace}

\title{\vspace{-3em} A Two-Stage Vehicle Routing Algorithm Applied to Disaster Relief Logistics after the 2015 Nepal Earthquake} 
\author{Stephanie Allen, \textit{SUNY Geneseo Mathematics Department \footnote{Contact information: stephanie.a.allen95@gmail.com}}\\
Advisor: Dr. Caroline Haddad, \textit{SUNY Geneseo Mathematics Department \footnote{Contact information: haddad@geneseo.edu}}}
\date{}

\begin{document}
\maketitle


\begin{abstract}
After the April 2015 Nepal Earthquake, the Himalayan Disaster Relief Volunteer Group distributed supplies to affected areas.  We model the organization’s delivery of supplies as a vehicle routing problem using Fisher and Jaikumar's two-stage method, which allocates locations to vehicles via an integer program and then uses heuristics to route the vehicles.  In the allocation stage, we use the assignment problem formulation to assign locations to vehicles.  In the routing stage, we implement multiple heuristics for the sake of comparison.  Our results illustrate the open nature of the vehicle routing problem and the computational necessity of heuristics.

\end{abstract}

\section{Introduction}

Operations research models can be used in disaster situations to improve efficiency as relief workers make critical decisions affecting the health and well-being of victims.  In this paper, we model the operations of the Himalayan Disaster Relief Volunteer Group (HDRVG), which formed after the April 2015 Nepal Earthquake and delivered supplies to affected areas for 26 days.  Using the group's detailed operational records, we seek to find the optimal assignment and routing decisions for the group's operations.  We model the logistical system of the organization as a Vehicle Routing Problem (VRP) because HDRVG's base of operations was ``a bed and breakfast'' called Yellow House, and it is reasonable to assume that vehicles departed from and returned to this location.  We choose to use Fisher and Jaikumar's two-stage VRP algorithm and provide background, theory, and results of the algorithm applied to HDRVG's operations. 
Section 2 examines potential models for disaster relief logistics.  Section 3 presents the details and theory of the methods used for the two-stage algorithm, drawing upon multiple sources. Section 4 discusses data, model assumptions, and model specification for the HDRVG network.  Section 5 discusses the results of the specified model for the HDRVG system, which reiterate the open nature of the VRP and the necessity of heuristics. Section 6 summarizes our findings and ideas for additional research.   

\section{Operations Research Models for Disaster Relief}

\subsection{Disaster Response Models}

Luis E de la Torre et al.'s literature review identifies some of the challenges in modeling disaster relief logistics which include: establishing the trade off between equity and speed, incorporating the uncertainty in crisis situations into models, incorporating the presence and/or number of command centers into models, incorporating types of aid, the chance of robbery, \& the loyalty/type of drivers into models, and using various types of objective functions \cite{luis2012disaster}. Caunhye et al.'s literature review looked at ``pre-disaster'' and ``short-term post-disaster'' logistic models.  In particular, for ``post disaster planning,” the models are ``mostly multi-period,” ``contemplate single objective[s],” ``rarely stochastic,” and ``construe relief distribution in terms of commodity flow or resource allocation,'' which means there are significant areas of open research for these models \cite{caunhye2012optimization}.  

Our model specifically fits into Caunhye et al.'s category of ``short-term post-disaster'' models as we seek to deliver supplies to locations right after the April 2015 Nepal Earthquake.  Both reviews demonstrate that there are plenty of open research areas in disaster relief modeling.


\subsection{General Network Models}

Before we review common transportation models, we present some network terminology for the purposes of this paper; all of the following definitions come directly from \cite{intro_graph_theory}.  A directed graph (digraph) is composed of nodes/vertices $V_{i}$ which are connected via arcs/edges $E_{i}$ that point in specified directions.  Each edge can be defined as an ordered pair of vertices $E_{i} = (V_{i},V_{j})$.  The network is composed of a vertex set $V$ and an edge set $E$.  More formally, the digraph has a vertex set $V$ and an irreflexive relation $E$ which is a subset of $VxV$.  A directed network is a function $G:E \rightarrow \mathbb{R}$, where the numeric result for our purposes represents the distance of the arc/edge. A path from node $v_{1}$ to $v_{n}$ in a graph $G$ is a sequence of nodes and edges such that no nodes are repeated. A cycle only repeats the beginning and ending nodes --- meaning $v_{1} = v_{n}$. Finally, a \textit{Hamiltonian Cycle} includes all the nodes of the graph $G$ exactly once except for the beginning node which is equal to the ending node.  A graph is \textit{Hamiltonian} if it has a \textit{Hamiltonian Cycle}.

Using this terminology, we define a general \textit{transportation model} as a set of ``supply'' and ``demand'' nodes where we seek to minimize the costs of transporting goods from the supply to the demand nodes \cite{LPOR_textbook}.  

The \textit{minimum cost flow model} accommodates intermediary nodes called transshipment nodes, which lie between the demand \& supply nodes and create more potential routes for goods \cite{LPOR_textbook}.

This is a powerful model because it can be modified to represent many transportation situations \cite{LPOR_textbook}. Finally, the \textit{vehicle routing problem model} represents situations in which we need to reach $n$ demand nodes using $m \geq 1$ vehicles such that the vehicles satisfy all of the demand nodes, do not overlap in their coverage of the nodes, and return to their mutual starting point (which we will refer to as the hub) \cite{intro_logistics_book,VRP_textbook,eksioglu2009vehicle}.  

Although all of these models can be used for disaster relief logistics, we choose the Vehicle Routing Problem Model because of the position of Yellow House as HDRVG's base of operations.

\section{Vehicle Routing Problem (VRP)}
The Vehicle Routing Problem (VRP) has been approached in many different ways, but no one has thus far developed a comprehensive solution technique that will guarantee an optimal solution for $n$ nodes \cite{VRP_textbook}.  
Therefore, we have to rely upon heuristics and metaheuristics.  

Within the VRP literature, the Traveling Salesman Problem (TSP) (which is the VRP but with $m=1$) has also been examined extensively, but no optimal algorithm has been devised.  Nevertheless, there are many heuristics available, which produce feasible solutions that differ in their distances from the true `optimal' solution.  These include modified minimum spanning tree algorithms such as Kruskal and Prim, and algorithms specifically for Hamiltonian circuits such as the Sub-tour Reversal algorithm, Nearest Neighbor Algorithm, and Sorted Edge Algorithm \cite{LPOR_textbook,MST,HamCycle,TSP_simple}.  General metaheuristics also exist that can be customized to the TSP, including the Tabu Search, Simulated Annealing, and Genetic Algorithms \cite{LPOR_textbook}.  There are also integer programming formulations of the TSP (some of which Pataki \cite{pataki2003teaching} illustrates in her article). Finally, Caric and Gold's book discusses more advanced formulations of the VRP \cite{VRP_textbook}.  

We use Fisher and Jaikumar's two-stage method to produce solutions for the formulation of the HDRVG network as a VRP (described in Brandimarte and Zotteri's \cite{intro_logistics_book}).  Fisher and Jaikumar's algorithm uses both integer programming techniques and heuristics/metaheuristics --- with the heuristics/metaheuristics drawing upon the TSP literature above.

\subsection{Two Stage VRP Algorithm: Overview}

In the two-stage method \cite{intro_logistics_book}, we ``assign'' nodes to vehicles using an integer program and then use heuristics/metaheuristics to find the optimal path for each vehicle to its assigned nodes.  For the assignment stage, we define $z_{ik}$ as the binary decision variable indicating whether or not node $i$ goes to vehicle $k$, $h_{ik}$ as the additional distance node $i$ would add to vehicle $k$'s route, $q_{i}$ as the amount to be delivered to node $i$, and $R_{k}$ as the capacity of vehicle $k$.  There are $m$ vehicles and $n$ demand nodes, so we have $nm$ decision variables.  Therefore, we can say:

\begin{equation}
\min Z = \sum\limits_{k=1}^{m} \sum\limits_{i=1}^{n} h_{ik} z_{ik}
\end{equation}

\begin{equation}
\sum\limits_{k=1}^{m} z_{ik} = 1\ \forall\ i=1...n 
\end{equation}

\begin{equation}
\sum\limits_{i=1}^{n} q_{i} z_{ik} \leq R_{k}\ \forall\ k = 1...m
\end{equation}

\begin{equation}
z_{ik} \in \{0,1\}
\end{equation}

\noindent The (2) constraints ensure coverage of each demand node by exactly one vehicle, and the (3) constraints enforce each vehicle's capacity.  The (4) constraints ensure the $z_{ik}$ variables are binary. The objective function (1) aims to minimize the $h_{ik}$ values. In order to calculate the $h_{ik}$ values --- the additional distance node $i$ would add to vehicle $k$'s route --- we need to use an approximation because, during the integer program solution process, the node assignments are unknown.  Therefore, we choose $m$ ``seed'' nodes which are ``far from the central location and far from each other'' such that $h_{ik}$ is the ``extra mileage'' of adding node $i$ to vehicle $k$'s route, with the route being represented by the line connecting the hub to the route's seed \cite{intro_logistics_book}.  The seeds serve as the other ends of a series of lines that come out from the hub.  Therefore we say $h_{ik} = d_{0i} + d_{i,\sigma_{k}} - d_{0,\sigma_{k}}$, with $d_{ij}$ as the distance between nodes $i$ and $j$ and with $\sigma_{k}$ as the $k$th seed.  Once the integer program produces assignments, we use heuristics and metaheuristics to route the vehicles \cite{intro_logistics_book}. Section 3.3 discusses our chosen heuristics/metaheuristics.

According to \cite{intro_logistics_book}, ``if there exists a feasible solution using $m$ vehicles [for the vehicle routing problem], we will find'' a feasible solution with the two-stage algorithm.  We guarantee this statement by proving its contrapositive.

\vspace{1em}

\noindent \underline{Theorem}: Assume a feasible solution with $m$ vehicles to the vehicle routing problem requires that the vehicles satisfy all $n$ demand nodes, do not overlap in their coverage of the nodes, and return to their mutual starting point \cite{intro_logistics_book,VRP_textbook,eksioglu2009vehicle}.  If the two-stage algorithm does not produce a feasible solution, there is no feasible solution for the vehicle routing problem with $m$ vehicles.

\vspace{1em}

\noindent \textit{Proof:} Given a set of nodes, we can, without loss of generality, write any ordering of the set of nodes, and this would be a feasible route.  Therefore, if the first stage of the algorithm produces feasible assignments of nodes, the second stage will produce a feasible route for each vehicle because, at a minimum, we could route the vehicles by choosing a random ordering of nodes such that Hamiltonian Cycles are formed.  Therefore, we must examine the first stage.

Consequently, if the two stage algorithm cannot produce a feasible solution, this means: 

\begin{equation}
\exists i = 1,...,n \text{ such that } \sum_{k=1}^{m} z_{ik} \neq 1
\end{equation}
or
\begin{equation}
\exists k = 1,...,m \text{ such that } \sum_{i=1}^{n} q_{i} z_{ik} > R_{k} 
\end{equation}

\noindent For the vehicle routing problem, we must be able to assign each node to one and only one vehicle to satisfy the feasibility requirements for a solution.  If (5) is true, then there is no feasible solution to the vehicle routing problem with $m$ vehicles because at least one node was not assigned to a vehicle or because at least one node was assigned to multiple vehicles.  We also know that the vehicle routing problem in a realistic setting requires that we satisfy the demand at each node.  However, if (6) is true, then our assignment has produced, for at least one vehicle, too much cargo such that it cannot deliver some of the supplies and thus cannot satisfy the demand of at least one node.  As a result, if (6) is true, there is no feasible solution to the vehicle routing problem given $m$ vehicles.  Consequently, we have shown that if we cannot find a feasible solution with the two-stage algorithm, we cannot find a feasible solution to the vehicle routing problem. QED

\vspace{1em}

As one last note for the overview of the two-stage method, for the ``extra mileage'' measure, the addition of another node into a route must either increase or have no effect on the length of the route \cite{intro_logistics_book}.  The Euclidean distance metric satisfies this condition because, for any three points $a$, $b$, and $c$ $\in \mathbb{R}^{2}$, $||a-b|| \leq ||a-c|| + ||c-b||$, as seen through the simple algebra of $||a-b|| = ||a-c+c-b|| \leq ||a-c|| + ||c-b||$ (due to the Euclidean distance triangle inequality) \cite{LA_textbook}.

\subsection{Stage One: Integer Programming for Node Assignment}

In the previous section, we explained the set-up of the integer program used in the first stage of the algorithm. This section will examine the method of solution used for this integer program. We use the \textit{MATLAB} mixed-integer linear programming solver, which uses a six step process to find a solution to a mixed-integer linear program; it uses the same procedure to solve a pure integer linear program such as ours.  We specify parameters that control these six steps.  

\subsubsection{Background Notation}

As we know, the decision variables in integer linear programs take on only integer values. Using the matrices $A$ and $C$, the column vectors $b$ and $d$, and the row vector $c$ with the knowledge that we have $n$ nodes and $m$ vehicles, we can write our integer program succinctly as \cite{IP_textbook}: 

\[\min \{ c \textbf{x}: \textbf{x} \in S\}, \hspace{1em} S = \{ \textbf{x} \in \{0,1\}^{nm}: A \textbf{x} = b \text{ and } C \textbf{x} \leq d \}\] 

\noindent The equality and inequality in $S$ are applied componentwise, and a solution is said to be feasible if and only if $\textbf{x} \in S$.  Because we have $n$ nodes and $m$ vehicles, we have $nm$ decision variables \cite{IP_textbook}.

\subsubsection{Six Step \textit{MATLAB} Integer Programming Solver Solution Process}

The information regarding these steps comes from the \textit{MATLAB} documentation \cite{matlab_doc_opts, matlab_doc_IP} as well as additional sources specified below.  According to the documentation for the \textit{intlinprog} solver, ``intlinprog can solve the problem in any of the stages. If it solves the problem in a stage, intlinprog does not execute the later stages'':

\begin{enumerate}

\item \textbf{Linear Program Pre-processing:} We choose to eliminate this step due to time constraints.  However, we do allow for pre-processing in Step 3.

\item \textbf{Relaxed Solution of Integer Program:} In this step, we eliminate the binary restriction on the variables (known as a relaxation of the integer program), and solve the problem as a linear program via the dual-simplex method \cite{LPOR_textbook, matlab_doc_IP}.  If the solution satisfies the binary constraints, we have our solution.  Otherwise, we continue to the next step.  

\item \textbf{Mixed-Integer Program Preprocessing:} ``Pre-processing'' steps ``try to tighten and simplify the formulation'' before utilizing other solution techniques \cite{IP_textbook}. If we have a constraint from our binary integer program, we can write it as \cite{IP_textbook}, 

\[ \sum_{j \in B^{-}} a_{j} x_{j} + \sum_{j \in B^{+}} a_{j} x_{j} \leq b_{i} \]

\noindent in which $B^{-} = \{j : a_{j} < 0\}$ and in which $B^{+} = \{j : a_{j} > 0 \} $.  We also can write the lower and upper bound of the left hand side of this constraint as,


\[ L_{min} = \sum_{j \in B^{-}} a_{j}, \hspace{1em} L_{max} = \sum_{j \in B^{+}} a_{j} \]

There are a few different preprocessing techniques (stated in \cite{IP_textbook}) which can be used that include:

\begin{itemize}

\item Infeasibility check: If $L_{min} > b$, the constraint isn't feasible

\item Redundancy: If $L_{max} \leq b$, we see the constraint could have a lower bound than $b$.

\item Variable Fixing: If there exists a $k \in B^{+}$ such that $L_{min} + a_{k} > b$, then $x_{k} = 0$.  If there exists a $k \in B^{-}$ such that $L_{min} - a_{k} > b$, then $x_{k} = 1$.

\item Improving Coefficients (using \cite{IP_textbook,savelsbergh1994preprocessing}): 

\begin{enumerate}

\item If $a_{k} > L_{max} - b$ for some $k \in B_{+}$, then we can subtract $(a_{k} - (L_{max} - b) )$ from $a_{k}$ and $b$.  We can do this because, if $x_{k} = 1$, then it won't matter if we subtract this quantity because we would just be subtracting some quantity from both sides of the inequality.  If $x_{k} = 0$, we know $a_{k} - (a_{k} - (L_{max} - b) )$ will be eliminated from the left hand side (LHS) but, since we would have on the right hand side (RHS):

\[ b - (a_{k} - (L_{max} - b)) = b - a_{k} + L_{max} - b = -a_{k} + L_{max} \]

\noindent we could add $a_{k}$ to both sides and leave $L_{max}$ on the RHS.  Therefore, the inequality would hold because the addition of $a_{k}$ to the LHS at most results in $LHS = L_{max}$.

\item If $a_{k} < b - L_{max}$, for some $k \in B_{-}$, we can replace $a_{k}$ with $b - L_{max}$.  If $x_{k} = 0$, the $a_{k}$ coefficient is irrelevant.  If $x_{k} = 1$, we can replace $a_{k}$ with $b - L_{max}$ because $a_{k} < b - L_{max}$ implies $a_{k} + L_{max} < b$.  Therefore, even if all variables with $k \in B_{+}$ are equal to 1, this $a_{k}$ coefficient alone will bring the LHS sum below the RHS $b$ value.  Consequently, replacing $a_{k}$ with $b - L_{max}$ would result in $a_{k} + L_{max} = b - L_{max} + L_{max} = b$ on the LHS, which our inequalities allow.

\end{enumerate}
\end{itemize}

\item \textbf{Cut Generation:} For ``cutting plane'' techniques, we solve the relaxed version of the integer program (IP): $\max \{ c \textbf{x}: \textbf{x} \in S\}$ subject to $S_{0} = \{ \textbf{x} \in \mathbb{R}^{nm}: A \textbf{x} = b \text{ and } C \textbf{x} \leq d \}$ \cite{IP_textbook}.  If the solution does not satisfy the binary constraints of $S$ ($S = \{ \textbf{x} \in \{0,1\}^{nm}: A \textbf{x} = b \text{ and } C \textbf{x} \leq d \}$), then we construct an additional constraint that the relaxed solution does not satisfy but that all $\textbf{x} \in S$ satisfy \cite{IP_textbook}.  This is a ``cutting plane'' because it ``cuts off'' the relaxed solution without removing any feasible binary solutions \cite{IP_textbook}. While many cutting plane algorithms exist, the \textit{MATLAB} solver utilizes the Gomory cutting plane algorithm for our integer programs, and it is as follows \cite{Gomory_bk}: 

\begin{enumerate}
\item Solve the relaxed version of the binary integer programming problem using the simplex method.  If a non-binary solution results, we proceed in the algorithm.

\item In the simplex tableau produced by the simplex method, look at the $i$th constraint $\sum_{j=1}^{n} t_{ij} x_{j} = x_{B_{i}}$ where $x_{B_{i}}$ is the value of the $i$th basic variable in the relaxed solution that is \textit{non-integer} and where the fractional part of $x_{B_{i}}$ is the largest.  Construct the cutting plane from this constraint via the following method:

	\begin{enumerate}
    \item Assuming $\lfloor d \rfloor$ is the ``integer part'' of $d$, we know $\lfloor t_{ij} \rfloor \leq t_{ij}$ such that $t_{ij} = \lfloor t_{ij} \rfloor + g_{ij}$ and know $\lfloor x_{B_{i}} \rfloor \leq x_{B_{i}}$ such that $x_{B_{i}} = \lfloor x_{B_{i}} \rfloor + f_{i}$ (with $g_{ij}$ and $f_{i}$ as the fractional parts).
    
    \item If a potential solution $\textbf{x}$ is an integer solution that satisfies the initial constraint, then it will satisfy a variation of this constraint: $\sum_{j=1}^{n} \lfloor t_{ij} \rfloor x_{j} \leq \lfloor x_{B_{i}} \rfloor$.
    
    \item We can use a slack variable to write $\sum_{j=1}^{n} \lfloor t_{ij} \rfloor x_{j} + u_{ij} = \lfloor x_{B_{i}} \rfloor$ where, if a potential solution $\textbf{x}$ is an integer, then $u_{ij}$ is an integer (indeed, we could also place an integer restriction upon $u_{ij}$).  Again, integer solutions that satisfied the initial constraint will satisfy this constraint.    
    
    \item If we subtract $\sum_{j=1}^{n} t_{ij} x_{j} = x_{B_{i}}$ from $\sum_{j=1}^{n} \lfloor t_{ij} \rfloor x_{j} + u_{ij} = \lfloor x_{B_{i}} \rfloor$, the result is $\sum_{j=1}^{n} -g_{ij} x_{j} + u_{i} = -f_{i}$ or, in other words, the Gomory cutting plane.  Since we subtracted an equal amount from both sides of the $\sum_{j=1}^{n} \lfloor t_{ij} \rfloor x_{j} + u_{ij} = \lfloor x_{B_{i}} \rfloor$ constraint, integer solutions that satisfied the initial constraint will satisfy this constraint.  Therefore, the constraint keeps the feasible integer solutions in the solution space.  
    
    \end{enumerate}

\end{enumerate}

\noindent To demonstrate how this new constraint eliminates a non-binary solution, we use \cite{Gomory_bk} and our knowledge of the simplex method.  We know \textit{basic variables} have coefficients of either 0 or 1 in the simplex tableau; therefore, when we construct the Gomory cutting plane, any non-zero $-g_{ij}$ coefficient values result from the original $t_{ij}$ values of the \textit{non-basic variables}, which are variables that are guaranteed to be 0 and which can have non-integer coefficients in the simplex tableau \cite{LPOR_textbook,Gomory_bk}.  The $-g_{ij}$ coefficient values for the basic variables are all 0 because the $t_{ij}$ coefficients are 0s or 1s.  Therefore, if we substitute the current solution into the Gomory cutting plane constraint, the result is $u_{i} = -f_{i}$.  However, we know $f_{i}$ must be positive because we assume that we found a feasible solution to the relaxed problem (which, as long as we assume the decision variables must be positive \cite{Gomory_bk}, would mean $x_{B_{i}}$ would have to be positive) and must be non-zero because we assume $x_{i}$ is not an integer (which would mean $x_{B_{i}} = \lfloor x_{B_{i}} \rfloor + f_{i}$, such that $f_{i} \neq 0$).  Therefore, since $u_{i} = -f_{i}$, then $\sum_{j=1}^{n} \lfloor t_{ij} \rfloor x_{j} > \lfloor x_{B_{i}} \rfloor$, which undermines the entire cutting plane inequality; furthermore, a slack variable should never be negative.  Consequently, the original solution does not satisfy this new constraint.

\item \textbf{Heuristics:} As the \textit{MATLAB} documentation states, ``there are techniques for finding feasible points faster before and/or during [the] branch-and-bound'' algorithm (to be mentioned below), and these are known as heuristic algorithms \cite{matlab_doc_IP}. We choose to use the Relaxation Induced Neighborhood Search (RINS) algorithm \cite{IP_textbook,danna2005exploring,matlab_doc_IP} which requires an initial feasible solution \cite{danna2005exploring,matlab_doc_IP}:  

\begin{enumerate}
\item We ``fix the variables that have the same values in the incumbent and in the current continuous relaxation'' \cite{danna2005exploring}.  The term ``incumbent'' refers to the best feasible solution.

\item We ``set an objective cutoff based on the objective value of the current incumbent'' \cite{danna2005exploring}.

\item We apply a combination of the branch-and-bound and cutting plane algorithms (see Steps 4 and 6) to the new formulation (known as the branch-and-cut algorithm).

\end{enumerate}

\item \textbf{Branch-and-Bound Algorithm:} For the branch-and-bound algorithm, we ``iteratively'' subdivide the above integer program into subproblems whose union is the entire set of feasible solutions and whose intersection is the empty set \cite{IP_textbook,LPOR_textbook}. Each time we divide the problem, we fix one of the variables such that it is set to $0$ in one subproblem and set to $1$ in the other subproblem.  Although the solver was free to use this technique, it never did as it solved the integer program for each of the 26 days.

\end{enumerate}

\noindent The table below summarizes the important parameter specifications made for the \textit{MATLAB} solver. 

\begin{table}[H]
\centering
\begin{tabular}{l|r}
\hline
\textbf{Parameter} & \textbf{Value} \\\hline
Cut Generation & Basic (solver only uses Gomory cuts for our IP) \\ Heuristics & Relaxation Induced Neighborhood Search (RINS) \\ IP Preprocess & Only a few \\ 
Root LP Algorithm & Dual-simplex \\ 
Tolerance Integer & 1e-5 (deviation from integer but still considered an integer) \\
\hline
\end{tabular}
\caption{\label{tab:widgets}Parameters for MATLAB Integer Linear Programming Solver}
\end{table}



\subsection{Stage Two: Meta-heuristics for Vehicle Routing}
Once the first stage is completed, we move into a series of Traveling Salesman Problems (TSP) because we must route each vehicle among its assigned nodes.  We employ multiple algorithms to route the vehicles and compare the results. 
First though, it is important to talk about the necessity of these algorithms when we could (in theory) check every permutation of the nodes through which a given vehicle could be routed.  This would not have been feasible for some of the routes.  For example, at maximum, it took 0.000068 seconds to calculate the length of a route and assign this length to a variable for Day 19 when the payload capacity was 2000kg.  Under a 2000kg payload capacity, Day 19's two routes consisted of 11 and 12 nodes.  Therefore, with this maximum time, to check all permutation for these routes, this would take:
\[ (11!/2) 0.000068 = 1357.17 \text{ seconds which is equivalent to 22.6 minutes} \]
\[ (12!/2) 0.000068 = 16286.1 \text{ seconds which is equivalent to 271.434 minutes} \]

\noindent We divide by two because, due to the circular natures of routes, an ordering of $1, 2, 3, 4, 5, 1$ would be equivalent to $1, 5, 4, 3, 2, 1$. Even more extreme, Day 10 had 14 nodes for one route, and a timing of 0.000062, which meant that it would take 31.279 \textit{days} to evaluate all of the permutations.  This computational time would be unacceptable in a disaster relief situation in which feasible solutions must be generated quickly. Therefore, as we can see from these calculations, we need heuristics and metaheuristics to save us time as we seek to find the best way of routing vehicles in disaster relief.

We first examine a simple \textit{Greedy Algorithm} whereby we build a vehicle route by adding nodes such that each addition minimizes the cumulative distance of the path up to that node.  The main disadvantage of this algorithm is that we do not consider the larger route when making our choices regarding which node to add into the route next.  For instance, the algorithm does not consider the fact that we will need to go back to the starting point to finish the route \cite{HamCycle,TSP_simple}. 

Next, we examine the \textit{Sub-Tour Reversal Algorithm}, which iteratively ``deletes exactly two links from [a] previous tour and replaces them by exactly two new links to form the new tour'' \cite{LPOR_textbook}.  This occurs by ``selecting a subsequence of the cities and simply reversing the order in which that subsequence of cities is visited'' \cite{LPOR_textbook}.  Given a sequence of nodes that form a Hamiltonian Cycle, $\{1,2,3,4,5,6,1\}$, some possible sub-tour reversals would be: $$\{1,3,2,4,5,6,1\}, \{1,2,3,6,5,4,1\}, \and \{1,5,4,3,2,6,1\}$$ 
\noindent The technique is well known and well documented \cite{LPOR_textbook,intro_logistics_book,other_subtour,usc_subtour,osman1993metastrategy}.  It is known as a ``local improvement procedure'' or a ``local search'' 
 algorithm because, although it improves a solution, it arrives at a local optimum rather than the global optimum \cite{LPOR_textbook,intro_logistics_book}.  In other words, the algorithm only focuses on its ``local neighborhood'' \cite{LPOR_textbook}.  The full algorithm proceeds as such \cite{LPOR_textbook}:

\begin{enumerate}

\item Establish an initial feasible solution.  For each vehicle, we choose the initial feasible solution as the nodes in numeric order with the Yellow House bed-and-breakfast as the first and last node (to form the Hamiltonian Cycle).

\item Create all possible sub-tour reversals using the current feasible solution by inverting all subsequences in the solution which will generate new routes (when compared to the current route).  Choose the new route that has the smallest distance as the new feasible solution.

\item Repeat Step 2 until sub-tour reversals do not result in smaller distances (when compared to the current feasible solution).

\end{enumerate}

\noindent Another disadvantage of this algorithm is the potential for computational issues as the number of nodes increases.  Given enough nodes, it may become too costly to check all feasible subsequences during multiple iterations. 

Finally, we look at the \textit{Simulated Annealing Algorithm} which, as opposed to the Sub-Tour Reversal algorithm, goes beyond just its local neighborhood to try to find the global optimum \cite{LPOR_textbook,osman1993metastrategy,usc_subtour}.  The algorithm follows the following procedure (from \cite{LPOR_textbook}):

\begin{enumerate}

\item Select a ``trial solution'' (which could also be termed a ``trial route'') whose objective function value is $z_{c}$.  For our application, our objective is to minimize the distance traveled by the vehicles, and the objective function value is the length of the route.  Also, the initial trial solution is the same as the initial trial solution for the Sub-Tour Reversal Algorithm.

\item  Utilizing the concept of a Sub-Tour Reversal, randomly choose the beginning and end of a subsequence to invert.  The beginning cannot be the first, the last, or the second to last node, the end cannot be the last node, and the subsequence itself cannot include both the second \textit{and} the second to last node.  All of these distinctions assume a route with the hub at the beginning and at the end of the route (or, in other words, the sequence of nodes).  

\item Find the length of the route generated from this Sub-Tour Reversal; call this length $z_{n}$.  

\begin{enumerate}

\item If $z_{n} \leq z_{c}$, then the route just generated from the Sub-Tour reversal becomes the new trial solution, and Step 2 repeats.

\item If $z_{n} > z_{c}$, the new route (generated from the Sub-Tour reversal) becomes the new trial solution if, given the parameter $T$ (called a ``temperature'') and a randomly chosen number $w$ from a uniform distribution with end points of 0 and 1, $e^{\frac{z_{c} - z{n}}{T}} > w$.  If this test fails, the original trial solution (before the Sub-Tour reversal) is the `new' trial solution, and we return to Step 2. 

\end{enumerate}

\item We continue Steps 2-3 based upon a chosen ``schedule'' of $T$s whereby we iterate a set number of times per $T$ \cite{LPOR_textbook}.

\end{enumerate}

\noindent We select large $T$ values at the beginning to generate a higher chance that we will utilize the new routes generated by the algorithm because ``the early emphasis is on taking steps in random directions...in order to explore as much of the feasible region as possible,'' which means the algorithm, with these larger $T$ values, will be more likely to ``accept'' solutions that have greater distances \cite{LPOR_textbook}.  The advantage to this strategy is that the algorithm will have a greater chance of ``escaping a local optimum'' \cite{LPOR_textbook}.  However, as time passes, we want to close in upon a ``good'' solution, so we decrease the $T$ values in order to prevent fewer ``worse'' solutions from being chosen \cite{LPOR_textbook}. We are not guaranteed to find the global optimum, but this algorithm increases our chances of ``escaping a local optimum.''



\section{Data, Model Assumptions, \& Model Specification}

\subsection{Data}

On April $25$, 2015 and May $12$, 2015, Nepal experienced two significant earthquakes and, throughout that time, experienced hundreds of after-shock quakes. The Himalayan Disaster Relief Volunteer Group (HDRVG) provided supplies to victims of the earthquakes via missions that brought food, water, shelter, and medicine to people across Nepal for 26 days \cite{FB_HDRVG,HDRN_data,HDRVG_website}.  HDRVG publicly released its data via the Tableau Public platform \cite{HDRN_data}.  The data set contains information regarding the date, location (District and Village Development Committee information including latitude and longitude information), mission number, type of product, and amount of product for each mission.  A mission can span multiple rows because a row represents one product that was delivered during a particular mission.

Using \textit{R} Version 3.2.0 with the dplyr and stringr packages, we wrote a script to extract information from this large data file and produced a CSV data file for each of the 26 days the organization ran missions.  Each of the data files contains the delivery VDC locations for that day.  We extracted the amount of kilograms of supplies distributed each day by writing a function to convert the unit types into a number of kilograms (such as sack = 25 kg), multiplying these conversions by the ``amount'' variable for each row, and summing the rows corresponding to each day to produce a file with the total supplies in kilograms available for each day.

\subsection{Model Assumptions}

We had to make some assumptions regarding our model in order to produce results during the time alloted to this research.  Therefore, we make the following assumptions/set the following parameters:

\begin{itemize}

\item We seek to minimize the distance traveled by the vehicles.
    
\item We assume all vehicles begin their journeys at the Yellow House Bed and Breakfast and return to Yellow House after making their routes.  These are reasonable assumptions because the organization refers to Yellow House as the central hub of activity for the group.
    
\item We treat each day independently, so we re-run the model each day using that day's data.  We do not consider any potential relationships between/among the circumstances of each day. 

\item For each day, we assume all supplies are available at the start of the day, that there is no difference between products, and that the products are distributed evenly across the locations.  These are not entirely realistic assumptions because a relief organization might be receiving supplies throughout the day, affected locations would need a variety of products (and thus would care about the difference between products), and some areas might need more goods than others (based on the degree to which different areas were affected).  Nevertheless, these assumptions enable us to formulate the problem.
    
\item We use the Euclidean Distance Metric to calculate the distances between nodes.  This does not take into consideration the road infrastructure of Nepal which would be important with trucks, but it does ensure that we do not violate the triangle inequality, which is important for our model.  In practice, our routes could be achieved with helicopters or drones.
    
\item We took the delivery locations established by the organization as given.  Therefore, we do not decide where to deliver supplies; we simply want to find the best way to route vehicles to the locations established by the organization.

\end{itemize}

\subsection{Model Specification}

We translated the two-stage method into code in \textit{MATLAB R2014b} with the Optimization Toolbox (Version 7.1) (which provided the integer program solver).  Our script executes the two-stage algorithm for all of the days of HDRVG's operations.  We had to select several parameter values for the two-stage method, which will be discussed in this section.  To begin with the integer program for the first stage, we discussed its formulation in Section 3.1 as follows: 

\[ \min Z = \sum\limits_{k=1}^{m} \sum\limits_{i=1}^{n} h_{ik} z_{ik} \]
\[ \sum\limits_{k=1}^{m} z_{ik} = 1\ \forall\ i=1...n \]
\[ \sum\limits_{i=1}^{n} q_{i} z_{ik} \leq R_{k}\ \forall\ k = 1...m \]
\[ z_{ik} \in \{0,1\} \]

\noindent First, we specified the $R_{k}$ value(s), which refer to the payload capacities of the vehicles (``the total weight of people and equipment...[a vehicle]...can carry without overloading'') \cite{payload_1}.  According to \cite{payload_2}, ``light-duty pickups'' have a payload capacity of about 900-1400 kilograms.  We chose to test two $R_{k}$ values, 1500 kg and 2000 kg, because the trucks in the HDRVG pictures appeared larger than ``light duty pickups'' \cite{FB_HDRVG}.  For the $q_{i}$ coefficients, we decided to assume each demand node in a given day would receive the same amount of supplies.  Consequently, for each day, we divided the total amount of supplies by the number of demand nodes to calculate the $q_{i}$ values for the day.  Finally, to determine the number of vehicles for a given day, we used the formula: $ceiling((n)/(payload/q))$. 

Next, as we discussed in Section 3.1, in order to calculate $h_{ik}$, we needed to use an approximation.  Therefore, we chose $m$ ``seed'' nodes which were ``far from the central location and far from each other'' such that $h_{ik}$ was the ``extra mileage'' of adding node $i$ to the route with the $k$th seed.  Consequently, we used the formula, $h_{ik} = d_{0i} + d_{i,\sigma_{k}} - d_{0,\sigma_{k}}$, where $d_{ij}$ was the distance between nodes $i$ and $j$ and $\sigma_{k}$ was the $k$th seed ($k = 1...m$) \cite{intro_logistics_book}.  The challenging part of this ``extra mileage'' formula was the determination of the ``seeds.''  We established a hierarchy of steps in the \textit{MATLAB} script to determine the seed locations for each day, with the number of seeds equal to the number of vehicles:

\begin{enumerate}
\item If the number of demand nodes is greater than 1, compute the convex hull of the nodes (including the hub) for the given day.  If the number of corners of the convex hull exceeds the number of vehicles $v$ for the day, choose $v$ corners.  

\begin{enumerate}

\item Otherwise, if the number of corners of the convex hull is less than the number of vehicles, take $v$ equally spaced points inclusively between the minimum and maximum latitude and longitude values of the entire set of nodes across all of the days.  

\end{enumerate}

\item If there is only one demand node, make that node the seed.
\end{enumerate}

\noindent Once we established the seed nodes for a given day, we calculated all of the $h_{ik}$ values for the objective function for that day.

For the heuristics/metaheuristics in stage two of the algorithm, only the simulated annealing algorithm had parameters to be specified.  The schedule of $T$ values --- with $z_{0}$ as the initial route distance --- was set according to \cite{LPOR_textbook} as: $T_{1} = 0.2z_{0}, T_{2} = 0.2T_{1}, T_{3} = 0.2T_{2}, T_{4} = 0.2T_{3}, \and T_{5} = 0.2T_{4}$.  We also ran the algorithm for 15 iterations per temperature $T_{i}$ for the 1500kg payload capacity and ran the algorithm for 20 iterations per temperature $T_{i}$ for the 2000kg payload capacity.


\section{Results}

For each of the 26 days, we produced six sets of results in order to determine the results using each of the heuristics/metaheuristics with both payload capacities (1500 kg and 2000 kg).  For the first stage, there is variability in the number of vehicles needed for each day.  Furthermore, some days' vehicles are assigned only one or two nodes while other days' vehicles are assigned more nodes.  The graphics below illustrate these observations by displaying the assignments of nodes to vehicles for a sample of days under each of the payload capacity situations (1500kg and 2000kg).  Nodes of the same symbol \textit{and} color belong to a common vehicle.

\begin{figure}[H]
\begin{center}
\includegraphics[height=0.33\textheight,keepaspectratio]{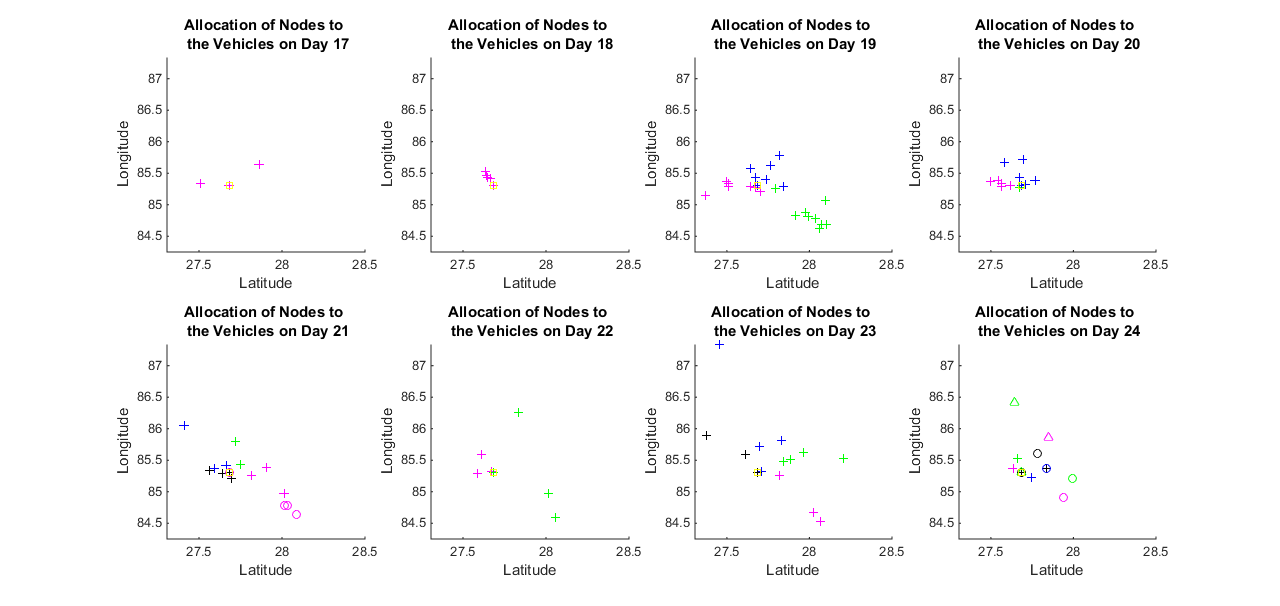}
\caption{\label{tab:widgets} Sample of Node Assignments for Eight Days with Payload = 1500kg}
\end{center}
\end{figure}

\begin{figure}[H]
\begin{center}
\includegraphics[height=0.33\textheight,keepaspectratio]{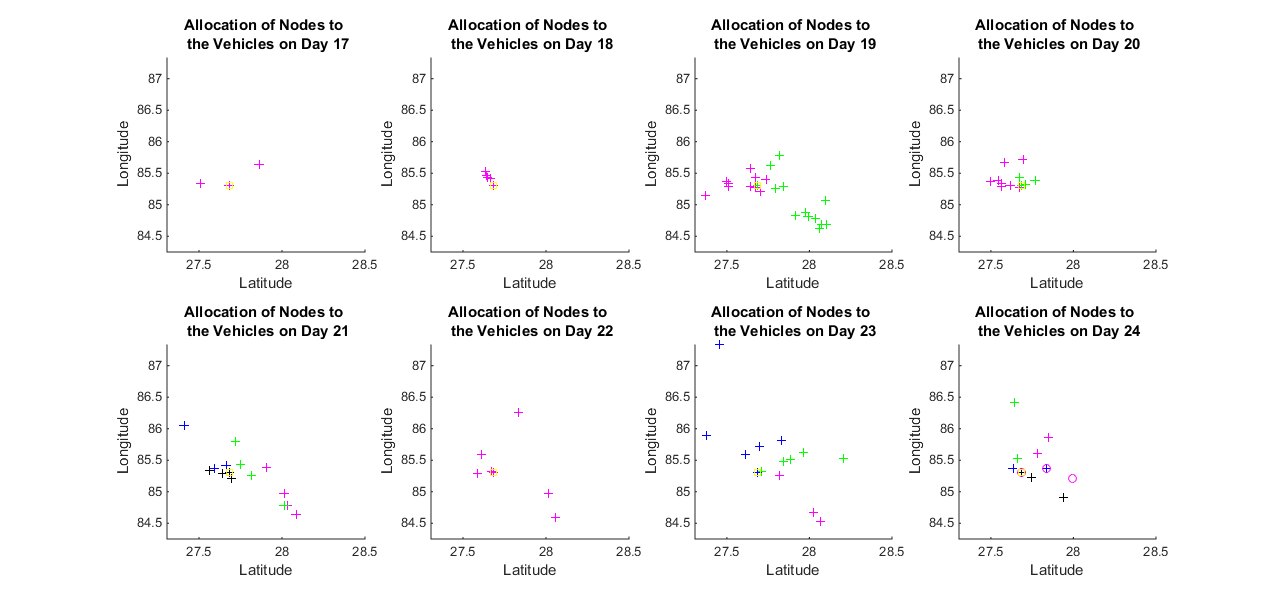}
\caption{\label{tab:widgets} Sample of Node Assignments for Eight Days with Payload = 2000kg}
\end{center}
\end{figure}

This sample illustrates the variability in assignments.  For Days 17 and 18, there are only a few nodes, and all are assigned to the same vehicle under both payload capacity values. 
For Days 19, 20, 21, and 23, there are quite a few nodes which cluster to certain vehicles.  There is also a difference in the number of vehicles and the resulting cluster patterns between the two payload capacities.  This is especially evident in Day 24; when the payload capacity of the vehicles is 1500kg, most of the vehicles are assigned only one node but, with the increase in the payload capacity to 2000kg, distinct clusters can be seen.  Indeed, for many of the days, it is evident that the number of vehicles decreases as the payload capacity increases.

For stage 2, to understand the differences between the three routing algorithms, it is useful to examine the routing of vehicles on Day 19 under each of the three algorithms for each of the two chosen payload capacities (which the graphs below illustrate).

\begin{figure}[H]
\begin{center}
\includegraphics[height=0.28\textheight,keepaspectratio]{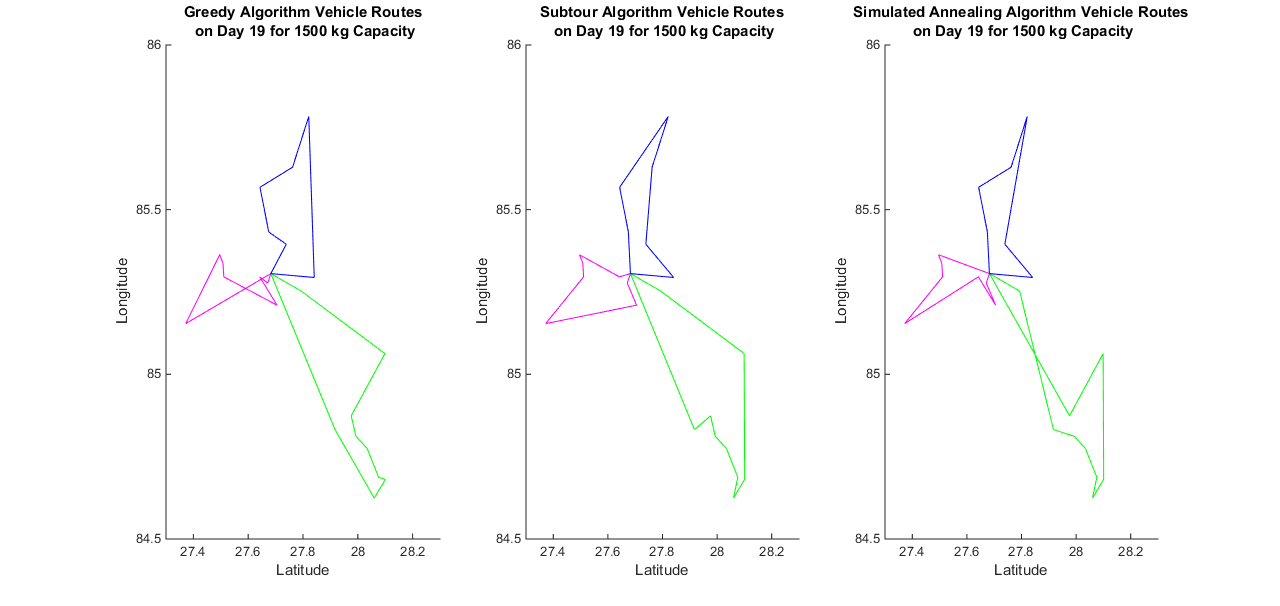}
\caption{\label{tab:widgets} Day 19 Routing under Three Algorithms with Payload = 1500kg}
\end{center}
\end{figure}

\begin{figure}[H]
\begin{center}
\includegraphics[height=0.28\textheight,keepaspectratio]{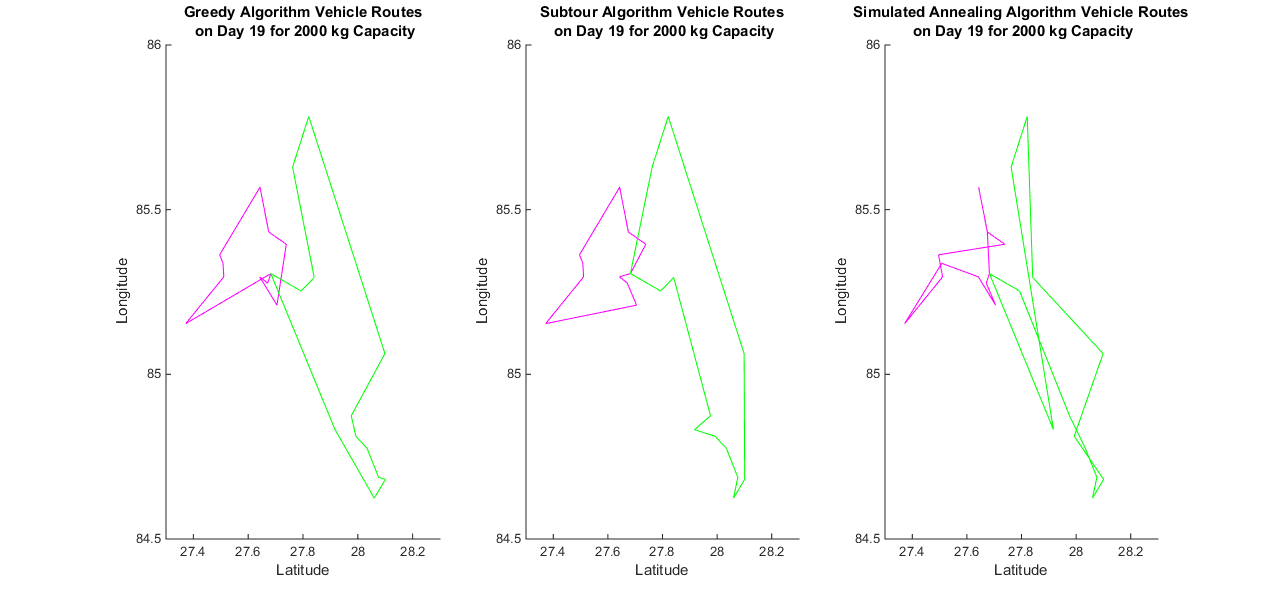}
\caption{\label{tab:widgets} Day 19 Routing under Three Algorithms with Payload = 2000kg}
\end{center}
\end{figure}

For each payload capacity, we see there are distinct differences in the shapes of the paths among each of the three algorithms. The crossings in some of the routes indicate inefficiencies, especially in the Simulated Annealing graph under the 2000kg payload capacity. 

To quantify our results, we found the total distance in degrees traveled by the vehicles each day under each heuristic/metaheuristic for each payload capacity.  We multiplied all of these degree measurements by 61.8 miles, which is the average of the degree-to-mile equivalences for latitude and longitude measures (1 degree latitude = 69 miles and 1 degree longitude = 54.6 miles) \cite{lat_long_miles}. We then summed these mile measurements to find the total number of miles traveled under each of the algorithms for each payload capacity.  The table below displays the results.

\begin{table}[H]
\centering
\begin{tabular}{l|c|c|c|c|c|c}
\hline
\textbf{Day} & \textbf{Greedy} & \textbf{Sub-Tour} & \textbf{Sim Annealing} & \textbf{Greedy} & \textbf{Sub-Tour} & \textbf{Sim Annealing} \\
   & \textbf{1500kg} & \textbf{1500kg} & \textbf{1500kg} & \textbf{2000kg} & \textbf{2000kg} & \textbf{2000kg} \\

\hline

Day 1 & 418.674 & 403.544 & 406.017 & 397.533 & 371.44 & 376.537 \\ Day 2 & 166.363 & 151.22 & 151.22 & 166.363 & 151.22 & 151.22 \\ Day 3 & 114.632 & 113.243 & 113.243 & 114.632 & 113.243 & 113.243 \\ Day 4 & 18.775 & 18.775 & 18.775 & 18.775 & 18.775 & 18.775 \\ Day 5 & 3.624 & 3.624 & 3.624 & 3.624 & 3.624 & 3.624 \\ Day 6 & 233.889 & 221.439 & 231.188 & 233.889 & 221.439 & 232.367 \\ Day 7 & 317.358 & 307.802 & 311.522 & 317.358 & 307.802 & 307.802 \\ Day 8 & 402.981 & 366.401 & 366.401 & 337.32 & 314.915 & 314.915 \\ Day 9 & 486.708 & 479.716 & 479.716 & 422.827 & 397.894 & 397.894 \\ Day 10 & 247.38 & 216.98 & 228.604 & 244.518 & 211.304 & 278.387 \\ Day 11 & 169.994 & 169.329 & 169.329 & 169.994 & 169.329 & 169.329 \\ Day 12 & 448.04 & 439.901 & 439.901 & 314.667 & 286.738 & 284.332 \\ Day 13 & 269.981 & 269.981 & 269.981 & 269.981 & 269.981 & 269.981 \\ Day 14 & 370.145 & 351.383 & 351.383 & 320.647 & 298.729 & 304.594 \\ Day 15 & 98.21 & 98.21 & 98.21 & 84.725 & 84.725 & 84.725 \\ Day 16 & 170.376 & 170.376 & 170.376 & 138.169 & 136.221 & 136.221 \\ Day 17 & 63.419 & 63.419 & 63.419 & 63.419 & 63.419 & 63.419 \\ Day 18 & 28.664 & 28.434 & 28.434 & 28.664 & 28.434 & 28.434 \\ Day 19 & 254.088 & 245.995 & 265.911 & 261.207 & 250.058 & 324.296 \\ Day 20 & 103.035 & 97.635 & 97.635 & 93.588 & 92.668 & 99.706 \\ Day 21 & 362.719 & 355.001 & 355.001 & 387.431 & 369.632 & 369.632 \\ Day 22 & 257.254 & 254.058 & 254.058 & 221.124 & 221.124 & 221.124 \\ Day 23 & 529.043 & 525.782 & 525.782 & 468.185 & 450.042 & 450.042 \\ Day 24 & 440.434 & 440.434 & 440.434 & 340.179 & 340.179 & 340.179 \\ Day 25 & 753.973 & 741.107 & 798.328 & 753.973 & 741.107 & 791.937 \\ Day 26 & 136.692 & 135.352 & 135.352 & 136.692 & 135.352 & 138.595 \\ \hline Totals & 6866.451 & 6669.142 & 6773.845 & 6309.484 & 6049.395 & 6271.31 \\

\hline
\end{tabular}
\caption{\label{tab:widgets} Total Distance Traveled by Vehicles Each Day Under the Three Routing Algorithms for 1500kg and 2000kg Payload Capacities}
\end{table}

For both payload capacities, the Sub-tour Reversal Algorithm is the best algorithm overall, followed by the Simulated Annealing Algorithm and then the Greedy Algorithm.  We expected the Simulated Annealing Algorithm to be the superior algorithm because it has a better chance of ``escaping a local minimum''; indeed, the algorithm is designed to be able to search a solution space \cite{LPOR_textbook}.  These results reconfirm that the Traveling Salesman Problem (TSP) has not been solved; as can be seen from the table, for any given day, there is not a guarantee regarding which algorithm will perform the best.  There are several days in the 2000 kg payload capacity chart in which the Greedy Algorithm produces a better solution (in terms of miles) than the Simulated Annealing Algorithm!

As we discussed in Section 3.3, it is not computationally feasible to examine all of the possible solutions when routing vehicles through several nodes.  These were the approximate run times for both stages of the algorithm under each of the heuristics and metaheuristics (for each of the capacities):

\begin{table}[H]
\centering
\begin{tabular}{l|c|c|c}
\hline
 & \textbf{Greedy} & \textbf{Sub-Tour} & \textbf{Sim Annealing} \\ \hline
 Run Time 1500kg (secs) & 1.650 & 2.066 & 2.670  \\
 Run Time 2000kg (secs) & 2.153 & 2.542 & 2.615 \\
 
\hline
\end{tabular}
\caption{\label{tab:widgets} Computation Time for Two Stage Algorithm for all 26 Days Under Each Routing Algorithm}
\end{table}

\noindent The speed of these algorithms when compared to the speed of checking permutations lends credibility to using these algorithms for vehicle routing problems.  In Section 3.3, we saw it takes 22.6 minutes to check 11 nodes, 271.434 minutes to check 12 nodes, and 31.279 days to check 14 nodes. 

If required to recommend one algorithm, we would recommend the Sub-Tour Reversal algorithm for the HDRVG's work because the Sub-Tour Reversal algorithm produces the best results under both payload capacities.  Our recommendation might change in the case of a significant increase in the number of nodes because the number of subsequences would increase with the number of nodes. 
With the Simulated Annealing Algorithm, we would be able to control the number of iterations directly, which would be useful in a potential situation in which the Sub-Tour Reversal Algorithm would be computationally infeasible.  However, given the computation time for these three algorithms, we would \textit{overall} recommend running all three algorithms because they would take 10 seconds or less to run (for a single payload capacity) and, as demonstrated by the charts above, one algorithm is not consistently superior for all of the days.

\section{Conclusion and Future Work}

This project demonstrates that operations research is a rich mathematical subfield filled with powerful theory and techniques that can allow us to model logistics scientifically.  Our goal was to minimize the distance traveled by vehicles as they brought supplies to affected communities.  We used a two-stage vehicle routing problem algorithm, which allowed us to explore integer programming techniques and heuristics/metaheuristics.  Our results reaffirmed the open nature of the vehicle routing problem and demonstrated that, until an optimal algorithm is found, it is necessary to use a variety of techniques and to then compare the output of these techniques to find the best solution.  

With regard to future work, there are quite a few complexities we would seek to integrate into the model.  We would seek to relax or even eliminate some of the assumptions present in this model, including the homogeneity of supplies \& distribution and the necessity of the Euclidean distance metric.  Instead of the Euclidean distance metric, we would aim to utilize road distances provided by Google Maps.  Integrating stochastic events into the model would also be extremely useful, especially given that the conditions right after natural disasters are often quite unclear.  Furthermore, we would research additional temperature scales for the Simulated Annealing algorithm to see if a different scale would improve the results of the algorithm for the HDRVG system. \footnote{Thanks to Joseph O'Brien for this suggestion at the 2017 Pi Mu Epsilon Conference.}

\bibliographystyle{plain}
\bibliography{main} 
\nocite{*}

\end{document}